


\input amstex
\documentstyle{amsppt}
\input cdlms.tex

\iftrue
\else
\baselineskip=1.67\normalbaselineskip     
\smallskipamount=1.33\smallskipamount     
\medskipamount=1.33\medskipamount         
\bigskipamount=1.33\bigskipamount
\magnification=\magstep1

\hsize=5.9truein
\hoffset=0.5truein
\catcode`@=11
\def\plainoutput{\shipout\vbox{\makeheadline\pagebody\vskip10pt\makefootline}%
  \advancepageno
  \ifnum\outputpenalty>-\@MM \else\dosupereject\fi}
\catcode`@=12
\fi

\iftrue
\magnification=\magstep1
\hsize=6.5truein
\hoffset=0.0truein
\baselineskip 1.4\normalbaselineskip

\tolerance=10000
\def\sqr{$\vcenter{\hrule height .3mm
\hbox {\vrule width .3mm height 2mm \kern 2mm
\vrule width .3mm} \hrule height .3mm}$}
\else

\baselineskip=1.33\normalbaselineskip
\tolerance=10000
\def\sqr{\vcenter{\hrule height .3mm
\hbox {\vrule width .3mm height 2mm \kern 2mm
\vrule width .3mm} \hrule height .3mm}}
\fi

\def \title{\medskip\centerline}
\def \Proof{\noindent {\bf Proof.\ \ }}
\def \R{{\Bbb R}}
\def \N{{\Bbb N}}

\def \a{\alpha}

\def \d{\delta}

\def \e{\epsilon}

\def \m{\mu}

\def \ind{\chi_k^m}

\font\ncbf=cmssbx10 at 9pt 
\def\abstract#1{{\ninepoint{\narrower\smallskip\noindent
	{\ncbf Abstract.} #1\smallskip}\vskip.5truein}}  
	
\def\hangbox to #1 #2{\vskip1pt\hangindent #1\noindent \hbox to #1{#2}$\!\!$}
 
\newskip\ttglue 
\def\ninepoint{\def\rm{\fam0\ninerm}
  	\textfont0=\ninerm \scriptfont0=\sixrm \scriptscriptfont0=\fiverm
  	\textfont1=\ninei  \scriptfont1=\sixi  \scriptscriptfont1=\fivei
  	\textfont2=\ninesy  \scriptfont2=\sixsy  \scriptscriptfont2=\fivesy
	\textfont3=\tenex  \scriptfont3=\tenex  \scriptscriptfont3=\tenex
	\textfont\itfam=\nineit  \def\it{\fam\itfam\nineit}
	\textfont\slfam=\ninesl  \def\sl{\fam\slfam\ninesl}
	\textfont\ttfam=\ninett  \def\tt{\fam\ttfam\ninett}
	\textfont\bffam=\ninebf  \scriptfont\bffam=\sixbf
	\scriptscriptfont\bffam=\fivebf  \def\bf{\fam\bffam\ninebf}
	\tt  \ttglue=.5em plus.25em minus.15em
	\normalbaselineskip=11pt
	\setbox\strutbox=\hbox{\vrule height8pt depth3pt width0pt}
	\let\sc=\sevenrm  \let\big=\ninebig \normalbaselines\rm}

	\font\ninerm=cmr9 \font\sixrm=cmr6 \font\fiverm=cmr5
	\font\ninei=cmmi9  \font\sixi=cmmi6   \font\fivei=cmmi5
	\font\ninesy=cmsy9  \font\sixsy=cmsy6 \font\fivesy=cmsy5
	\font\nineit=cmti9  \font\ninesl=cmsl9  \font\ninett=cmtt9
	\font\ninebf=cmbx9  \font\sixbf=cmbx6 \font\fivebf=cmbx5
	\def\ninebig#1{{\hbox{$\textfont0=\tenrm\textfont2=\tensy
	\left#1\vbox to7.25pt{}\right.$}}}

\font\efont=cmti8
\font\adfont=cmr8

\title {\bf A Lifting Theorem for Locally Convex Subspaces of $L_0$}
\smallskip
\author R. G. Faber
\endauthor
\centerline{R. G. Faber}

\bigskip

\address 
\ 
\line{Department of Mathematics \hfil} 
\line{University of Illinois \hfil}
\line{Urbana, IL 61801 \hfil}
\line{U.S.A. \hfil}
\line{\efont E-mail: \adfont faber\@math.uiuc.edu \hfil}
\endaddress

\abstract {We prove that for every closed locally convex subspace $E$ of $L_0$ and
for any continuous linear operator $T$ from $L_0$ to $L_0/E$ there is a
continuous linear operator $S$ from
$L_0$ to $L_0$ such that $T = QS$ where $Q$ is the quotient map from
$L_0$ to $L_0/E$.}

\footnote""{1991 {\it Mathematics Subject Classification}: Primary 46E30, 46A22.}
\footnote""{This paper is part of the author's Ph.D. thesis at the University of
Illinois, written under the supervision of N.T. Peck. The research was supported 
in part by a Beckman grant from the Research
Board of the University of Illinois.}

\noindent {\bf \S 0. Introduction.}

Let $E$ be a subspace of $L_0 = L_0[0,1]$, the space of all measurable
functions from $[0,1]$ to $\R$. Let $T$ be an operator from $L_0$ to
$L_0/E$.  What conditions on $E$ ensure that we can find an operator $S$
that makes the following diagram commute?
$$
\Cgaps{1.6}
\CD
                                       &     L_0 @()\l{Q} @(0,-1) \\
L_0 @()\L{T} @(1,0) @() \a-\L{S} @(1,1) 
   &     L_0/E
\endCD
$$
A. Pe\l czy\'nski was the first to ask if locally convex subspaces $E$ have this
property.
If $E$ is locally bounded then we can find such an operator (Kalton - Peck
\cite{2}).
Peck - Starbird \cite{6} showed that this is also true when $E$ is isomorphic
to $\omega$, the space of all real sequences.  The goal of this paper
is to show that if $E$ is locally convex then we can complete the
previous diagram. 

We will state some notation. We will let $\mu$ represent the standard
Lebesgue measure.  We also define the map $f \mapsto \|f\|_0 \ \ (L_0 \to
\R)$ as 
$$
\|f\|_0 = \int_0^1{|f(x)| \over 1 + |f(x)|}\,dx.
$$
This map is an F-norm on $L_0$, that is,
$$
\matrix
\hbox{\rm (i)} \hfill & \|f\|_0 > 0 \hfill & \hbox{\ \ \ \ \ \ \ } & f \neq 0,
 \hfill & \hbox{\ \ \ \ \ \ \ \ \ \ \ \ \ \ \ \ \ \ \ \ \ } \\
\hbox{\rm (ii)} \hfill & \|\alpha f\|_0 \leq \|f\|_0 \hfill & & |\alpha|\leq 1
 \hbox{\rm \ and } f\in L_0, \hfill & \\
\hbox{\rm (iii)} \hfill & \lim_{\alpha \to 0} \|\alpha f\|_0 = 0 \hfill & & f \in L_0, 
 \hfill & \\
\hbox{\rm (iv)} \hfill & \|f + g\|_0 \leq \|f\|_0 + \|g\|_0 \hfill & &
 f,g \in L_0. \hfill & \\
\endmatrix
$$

\noindent
The map also induces a metric on $L_0$.  The topology induced by the
$L_0$ metric is just the topology of convergence in measure.
For $f \in L_0$ we define $\sigma: L_0 \to [0,1]$ by
$$\sigma(f) = \sup_{n \in \N} \|nf\|_0 = \lim_{n \to \infty} \|nf\|_0.
$$
By the dominated convergence theorem we can see that $\sigma(f) =
\m(\hbox{\rm supp\ }f)$, where supp~$f = \{x:|f(x)| > 0\}$.
The F-norm on the quotient space $L_0/E$ is defined in the usual way

$$
\|\gamma\|_{L_0/E} = \inf_{f \in \gamma} \|f\|_0 
\hbox{\rm\ \ \ \ \ \ for all\ } \gamma \in L_0/E.
$$
For a subset $A$ of $[0,1]$ we will let $L_0(A)$
mean the subspace of $L_0$ consisting of all functions
supported on $A$. We define
$$
\|f\|_{L_0(A)} = \|f\cdot \chi_A\|_0,
$$
where $\chi_A$ is the characteristic function of $A$.
\bigskip

\noindent{\bf \S 1. Preliminary Lemmas.}

We have a lifting theorem for locally bounded subspaces (See Theorem 3.6
of 
\cite{2}.) and we will see
that locally convex subspaces are in some sense almost locally bounded.  
The lemmas that follow show us that the `unbounded part' of a locally
convex subspace is arbitrarily small. Lemma 1.2 is at the heart of this
argument. However, we first need a lemma from Paley and Zygmund \cite{5}.
\bigskip


\proclaim {Lemma 1.1} Let $\alpha > \beta \geq 0$. If $f \in L_0[0,1]$
such that $\int_0^1 f \geq \alpha$ and $\|f\|_2 = 1$ then $$\mu
(t:f(t) \geq \beta ) \geq (\alpha - \beta)^2.$$
\endproclaim

\Proof 
$$
\eqalignno{
\alpha &\leq \int_0^1f = \int_{\{f\geq \beta\}}f +
\int_{\{f<\beta\}}f \cr
&\leq \int_0^1 f\cdot I_{\{f\geq \beta\}} + \beta \cr
&\leq
\|f\|_2\cdot \|I_{\{f\geq\beta\}}\|_2 + \beta & (\hbox{\rm Schwarz Inequality})
\cr
&= \sqrt{\mu(t:f(t)
\geq\beta)} + \beta.
}$$
\line{{\rm So } $\mu(t:f(t)\geq \beta) \geq (\alpha - \beta)^2.$
 \hfil $\square$}
\bigskip

Notice that Rademacher functions do not appear in the statement of the
next lemma but they play a key role in the proof.  Recall that all the
Rademacher functions act on $[0,1]$ and have values in the two point set
\{-1,1\}. The first Rademacher function, $r_1$, is $1$ everywhere. The
second, $r_2$, is $1$ on $[0,{1 \over 2})$ and $-1$ on $[{1\over2},1]$; $r_3$
is $1$ on $[0,{1\over4})$ and $[{1\over2},{3\over4})$
but $-1$ on $[{1\over4},{1\over2})$
and $[{3\over4},1]$; and so on.
For convenience we will say that a sequence of functions
$(f_i)_{i=1}^\infty$ is {\it $\d$-tapering} if $\|2^i \cdot f_i\|_0
\leq \d$ for all $i \geq 1$.
\bigskip


\proclaim {Lemma 1.2} Let $E$ be a locally convex subspace of $L_0$. 
For every $\epsilon > 0$ there is a $\delta > 0$
such that if $(f_i)_{i=1}^\infty \subset E$ is $\d$-tapering
then
$$\mu \left(\, \bigcup_{i=1}^\infty\left\{x:|f_i(x)| > 1\right\}\right)
\leq \epsilon.$$
Moreover, $\d$ can be chosen to be any positive number such that the
closed convex hull of $\{f\in E:\|f\|_0 \leq \d\}$ is contained in 
$\{f\in L_0:\|f\|_0 \leq {\e \over 80}\}$. 
\endproclaim

\Proof Consider the following function on $[0,1]$:
$$
g(t) = \left|{1 \over \sqrt{\sum_{j=1}^Na_j^2}}
\sum_{k=1}^Na_kr_k(t)\right|,
$$
where $a_1,a_2,\cdots,a_N \in \R$ and $r_1,r_2,\cdots,r_N$ are the
first $N$ Rademacher functions. Then from Khinchine's inequality we
have
$$\eqalign{
\int_0^1 g &= {1 \over \sqrt{\sum_{j=1}^N a_j^2}} \cdot \left\|\sum_{k=1}^N
a_k\cdot r_k\right\|_1 \cr
&\geq {1 \over \sqrt{\sum_{j=1}^N a_j^2}} \cdot {1 \over
2}\sqrt{\sum_{k=1}^Na_k^2} \cr
&= {1 \over 2}.
}$$
Since the Rademacher functions are orthonormal over $[0,1]$ we have
$$\eqalign{
\|g\|_2^2 &= \left\|\sum_{k=1}^N \left({a_k \over
\sqrt{\sum_{j=1}^Na_j^2}}r_k\right) \right\|_2^2 \cr
&= \sum_{k=1}^N \left({a_k \over \sqrt{\sum_{j=1}^Na_j^2}}\right)^2 \cr
&= 1.
}$$
We now are ready to use Lemma 1.1 with $\alpha = 1/2$ and $\beta = 1/4$:
$$\mu\left(t:{1 \over \sqrt{\sum_{j=1}^N a_j^2}}
  \left|\sum_{k=1}^N a_k r_k(t)\right| \geq {1 \over 4}\right) \geq 
\left({1 \over 2} -
{1 \over 4}\right)^2 = {1 \over 16}.$$
Therefore,
$$\mu\left(t:\left|\sum_{k=1}^N a_k r_k(t)\right| \geq 
  {1 \over 4} \sqrt{\sum_{j=1}^N a_j^2} \,
\right) \geq {1 \over 16}\eqno (*)$$
for $a_1,a_2,\cdots,a_N \in \R$.

Let $\epsilon > 0$ be given. Since $E$ is locally convex there is a 
$\delta > 0$ such that the closed convex hull of $\{f\in E:\|f\|_0 \leq
\delta\}$ is contained in $\{f \in L_0:\|f\|_0 \leq {\epsilon \over
80}\}$. Suppose $(f_i)_{i=1}^\infty \subset E$ is $\d$-tapering.
Then for every $N \geq 0$ we have
$$\eqalignno{
\epsilon/80 &= \int_0^1 \epsilon/80 \, dt \cr
&\geq \int_0^1 \left\|\sum_{i=1}^N {1 \over 2^i} r_i(t) 2^i f_i(x)
\right\|_0\,dt
& \hbox{\rm (local convexity)} \cr
&= \int_0^1 \int_0^1 {|\sum_{i=1}^Nr_i(t)f_i(x)| \over 
1 + |\sum_{i=1}^Nr_i(t)f_i(x)|}\,dx\,dt \cr 
&= \int_0^1 \int_0^1 {|\sum_{i=1}^Nf_i(x)r_i(t)| \over
1 + |\sum_{i=1}^Nf_i(x)r_i(t)|}\,dt\,dx & \hbox{\rm (Tonelli)} \cr
&\geq {1 \over 16}\int_0^1{{1 \over 4} \sqrt{\sum_{i=1}^Nf_i(x)^2}
\over 1 + {1 \over 4}\sqrt{\sum_{i=1}^Nf_i(x)^2}}\,dx.
 & \hbox{\rm (by $*$)} \cr
}$$ 
So 
$$
\int_0^1{{1 \over 4} \sqrt{\sum_{i=1}^Nf_i(x)^2}
\over 1 + {1 \over 4}\sqrt{\sum_{i=1}^Nf_i(x)^2}}\,dx \leq 16 {\epsilon
\over 80} = {\epsilon \over 5}.
$$
for all $N \geq 1$.
Therefore $\mu(x:\sqrt{\sum_{i=1}^Nf_i(x)^2} > 1) \leq \epsilon$
for all $N$. Indeed, suppose not; then
$$
\int_0^1{{1 \over 4} \sqrt{\sum_{i=1}^Nf_i(x)^2} 
\over 1 + {1 \over 4}\sqrt{\sum_{i=1}^Nf_i(x)^2}}\,dx > \epsilon \, 
\left({1/4
\over 1 + 1/4}\right)
= {\epsilon \over 5}.
$$
This is a contradiction.
Thus $\mu (x: \sum_{i=1}^N f_i(x)^2 > 1) \leq
\epsilon$ for all $N \geq 1$.
Letting $N$ go to infinity we get $\mu(x:\sum_{i=1}^\infty f_i(x)^2
> 1) \leq \epsilon$. Finally, since 
$$
\bigcup_{i=1}^\infty \{x:|f_i(x)| > 1\} \subset
\left\{x:\sum_{i=1}^\infty f_i(x)^2 > 1\right\}
$$
we can conclude that
$$
\mu\left(\,\bigcup_{i=1}^\infty\left\{x:|f_i(x)| > 1\right\}\right) 
\leq \epsilon. \eqno
\square
$$ 
\bigskip

Lemmas 1.3 and 1.4 find an arbitrarily small set that contains all the
`unboundedness' of $E$.

\bigskip


\proclaim{Lemma 1.3} Let $E$ be a locally convex subspace of $L_0$. Let 
$\e > 0$ and find $\d > 0$ so that the closed convex hull of 
$\{f\in E:\|f\|_0 \leq \d\}$ is contained in $\{f\in L_0: \|f\|_0\leq
{\e \over 80}\}$. Then if $((f_i^{(k)})_{i=1}^\infty)_{k=1}^\infty
\subset E$
is any countable
collection of $\d$-tapering
sequences then
$$
\mu\left(\,\bigcup_{k=1}^\infty \bigcap_{l=1}^\infty
\bigcup_{i=l}^\infty \{x:|f_i^{(k)}(x)| > 1 \}\right) \leq \e.
$$
\endproclaim

\Proof
We will start
by considering the first $n$ sequences. Let $N_1 < N_2 < \cdots
< N_{n-1}.$ Then
$$
(f_i^{(1)})_{i=1}^{N_1} \cup (f_i^{(2)})_{i=N_1+1}^{N_2} \cup
\cdots \cup (f_i^{(n-1)})_{i=N_{n-2}+1}^{N_{n-1}} \cup
(f_i^{(n)})_{i=N_{n-1} +1}^\infty
$$
is another $\d$-tapering sequence.
So for all $N_{n-1}$ we have
$$
\eqalign{
\e &\geq \m\biggl(\,\bigcup_{i=1}^{N_1}\{x:|f_i^{(1)}(x)| > 1\} \cup
\bigcup_{i=N_1 + 1}^{N_2}\{x:|f_i^{(2)}(x)| > 1\} \cup 
\cdots
\cr
& \ \ \ \ \ \ \ \ \ \ \ \ \ \ \ \ 
\cup
\bigcup_{i=N_{n-2} + 1}^{N_{n-1}}\{x:|f_i^{(n-1)}(x)| > 1\} \cup
\bigcup_{i=N_{n-1} + 1}^\infty \{x:|f_i^{(n)}(x)| > 1\}
\biggr) \cr
&\geq
\m\biggl(\,\bigcup_{i=1}^{N_1}\{x:|f_i^{(1)}(x)| > 1\} \cup 
\bigcup_{i=N_1 + 1}^{N_2}\{x:|f_i^{(2)}(x)| > 1\} \cup 
\cdots 
\cr
& \ \ \ \ \ \ \ \ \ \ \ \ \ \ \ \ 
\cup 
\bigcup_{i=N_{n-2} + 1}^{N_{n-1}}\{x:|f_i^{(n-1)}(x)| > 1\} \cup 
\bigcap_{l=1}^\infty \bigcup_{i=l}^\infty \{x:|f_i^{(n)}(x)|
 > 1\} 
\biggr). 
}$$
Let $N_{n-1}$ go to infinity to obtain 
$$\eqalign{
&\m\biggl(\,\bigcup_{i=1}^{N_1}\{x:|f_i^{(1)}(x)| > 1\} \cup
\bigcup_{i=N_1 + 1}^{N_2}\{x:|f_i^{(2)}(x)| > 1\} \cup
\cdots
\cr
& \ \ \ \ \ \ \ \ \ \ \ \ \ \ \ \ 
\cup
\bigcup_{i=N_{n-2} + 1}^\infty\{x:|f_i^{(n-1)}(x)| > 1\} \cup
\bigcap_{l=1}^\infty \bigcup_{i=l}^\infty \{x:|f_i^{(n)}(x)|
 > 1\}
\biggr) \leq \e.
}$$
Repeat this step $n-2$ times to get
$$
\m\left(\,\bigcup_{k=1}^n \bigcap_{l=1}^\infty 
\bigcup_{i=l}^\infty \{x:|f_i^{(k)}(x)| > 1\} \right) \leq \e.
$$
Let $n$ go to infinity to get the desired conclusion,
$$
\m\left(\,\bigcup_{k=1}^\infty \bigcap_{l=1}^\infty
\bigcup_{i=l}^\infty \{x:|f_i^{(k)}(x)| > 1\} \right) \leq \e. \eqno 
\square
$$

\bigskip
 
In the proof of Lemma 1.4 we use the fact that the space of all
Lebesgue measurable subsets of $[0,1]$ is a complete separable metric 
space. The distance definition is 
$$d(A,B) = \m(A \vartriangle B),$$
where $A \vartriangle B$ stands for the symmetric difference $(A~\setminus~B)
\cup (B~\setminus~A)$. 
We consider $A$ and $B$ to be identical if $\m(A \vartriangle B) = 0$.

\bigskip


\proclaim{Lemma 1.4} Let $E$ be a locally convex subspace of $L_0$. Let 
$\e > 0$ and find $\d > 0$ such that the closed convex hull of 
$\{f\in E: \|f\|_0 \leq \d\}$ is contained in 
$\{f\in L_0: \|f\|_0 \leq {\e \over 80} \}$.  Then there is a measurable
set $A$, $\m(A) \leq \e$, such that if $(f_i)_{i=1}^\infty \subset E$ is any 
$\d$-tapering sequence then
$$
\m \left(\,\bigcap_{l=1}^\infty \bigcup_{i=l}^\infty
\{x:|f_i(x)| > 1\} \setminus A\right) = 0.
$$  
\endproclaim
 
\Proof
Let $(f_i^{(t)})_{i=1}^\infty, t \in T$, be the collection of all 
sequences in $E$ such that $\|2^i \cdot f_i^{(t)}\|_0 \leq \d$ 
for all $i$.  $T$ could be an uncountable index set.  For each 
$t \in T$ define
$$
A_t = \bigcap_{l=1}^\infty \bigcup_{i=l}^\infty
\{x:|f_i^{(t)}(x)| > 1 \}.
$$
$(A_t)_{t \in T}$ is a subspace of the separable metric space 
consisting of all Lebesgue measurable subsets of $[0,1]$.  So $(A_t)_{t\in T}$
is separable.  Let $(A_{t_j})_{j=1}^\infty$ be a countable dense
subset.  Let 
$$
A = \bigcup_{j=1}^\infty A_{t_j}.
$$
By Lemma 1.3 $\m(A) \leq \e$. 
Let $\eta > 0$ and $t \in T$ be given. There is a $j$ 
such that $\m(A_{t_j}~\vartriangle~A_t)~<~\eta$ since 
$(A_{t_j})_{j=1}^\infty$ is dense in $(A_t)_{t \in T}$.
$$
\m(A_t~\setminus~A) \leq \m(A_t~\setminus~A_{t_j})
\leq \m(A_t~\vartriangle~A_{t_j}) < \eta.
$$
\line{Since $\eta > 0$ is arbitrary, $\m(A_t~\setminus~A) = 0$ 
for all $t \in T$. \hfil $\square$}

\bigskip

We are now ready to prove the main theorem.  The proof for locally 
bounded spaces in Kalton, Peck and Roberts \cite{3} was the inspiration
for this proof.  However, the proofs are quite different in places.
\bigskip

\vfil
\eject


\noindent{\bf \S 2. The Lifting Theorem.}

\proclaim{Theorem 2.1} Let $E$ be a locally convex subspace of $L_0[0,1]$.
Let $T:L_0[0,1] \to L_0[0,1]/E$ be a continuous linear operator.  
Then there is a unique continuous linear operator $S:L_0[0,1] \to L_0[0,1]$ so that
$T=QS$, where $Q:L_0[0,1] \to L_0[0,1]/E$ is the quotient map.
$$
\Cgaps{.8}
\CD
                                      &  &     L_0 @()\l{Q} @(0,-1) \\
L_0 @()\L{T} @(2,0) @()\a-\L{S} @(2,1) 
  &  &     L_0/E
\endCD
$$
\endproclaim

\Proof 
For each $n=1,2,3,\cdots$ 
find $\d_n > 0$ so that the closed convex hull of $\{f\in E:\|f\|_0 \leq \d_n\}$
is contained in $\{f\in L_0:\|f\|_0 \leq 1/80n\}$, and use Lemma 1.4 to find a
measurable set $A_n$ so that
$$
\eqalign{
\ \ &(i)\  \m(A_n) \leq {1 \over n}, \cr
&(ii)\ (f_i)_{i=1}^\infty \subset E \hbox{\rm\ and\ }
\left\|2^i\cdot f_i\right\|_0 \leq \d_n \hbox{\rm\ for all\ }i 
\Rightarrow \m\left(\, \bigcap_{l=1}^\infty \bigcup_{i=l}^\infty
\{x:|f_i(x)| > 1\} \setminus A_n \right) = 0.
}$$
Without loss of generality we may assume that 
$\d_1 \geq \d_2 \geq \d_3 \geq \cdots\!,\ \d_n \to 0$, and referring to the 
construction we can take
$A_1 \supset A_2
\supset A_3 \supset \cdots$.
Since $T$ is continuous, for each $\d_n$ we can find $\e_n>0$ so
that $\|f\|_0 \leq \e_n \Rightarrow \|Tf\|_{L_0/E} \leq \d_n / 6$.
Without loss of generality we may also assume that $\e_1 \geq \e_2
\geq \e_3 \geq \cdots$.
For each $m$ and 
$k=1,2,\cdots 2^m$ define 
$$
\Delta_k^m = \left[{k-1 \over 2^m},{k \over
2^m}\right).
$$
Define 
$$
\chi_k^m = \chi_{\Delta_k^m} \hbox{\rm \ for } k=1,\cdots,2^m;\
m=1,2,\cdots.
$$

Let $v \in L_0$ be given.  Define $S(0) = 0$.  So we will assume $v \neq 0$.
For the next few pages we will work
to define $S(v)$. Define $w_k^m = v\cdot \chi_k^m$ for $k=1,\cdots,2^m;
\ m=1,2,\cdots.$
For the time being we will consider $m$ and $k$ to be fixed and look
at $w_k^m$.
Let $m_0$ be the smallest integer so that $1/2^{m_0} \leq \e_1$, 
and assume $m \geq m_0$. Let $n(m)$ be the largest integer so that
$\e_{n(m)} \geq 1/2^m$. Since $T$ is continuous we know that $n(m)$ goes
to infinity as $m$ goes to infinity unless $T$ is identically $0$. 
For each $i=1,2,\cdots$ we can select $g_i \in L_0$ so that
$Qg_i = Tw_k^m$ and
$$
\left\|4^i\cdot g_i\right\|_0 \leq \left(1 + {1 \over 2i}\right)
\left\|4^i \cdot Tw_k^m\right\|_{L_0/E}.
$$
If $v=0$ then $g_i = 0$ for all $i = 1,2,\cdots$.
Note that $\|4^i\cdot w_k^m\|_0 \leq 1/2^m \leq \e_{n(m)}$ for all
$i=1,2,\cdots$. So $\|4^i\cdot Tw_k^m\|_0 \leq \d_{n(m)}/6$ for all
$i=1,2,\cdots$. Therefore $\sigma(Tw_k^m) \leq \d_{n(m)}/6$.
For $j \geq i \geq 1$,
$$
\left\|4^i(g_i - g_j)\right\|_0 \leq \left\|4^i\cdot g_i\right\|_0
+ \left\|4^j\cdot g_j\right\|_0 \leq \left(2 + {1\over 2i} 
+{1 \over 2j}\right)\cdot \sigma\left(Tw_k^m\right) \leq
3 \cdot {\d_{n(m)} \over 6} < \d_{n(m)}.
$$
Let $f_i = 2^i(g_i - g_{i+1})$.
Then $\|2^i \cdot f_i\|_0 \leq \d_{n(m)}$ for all $i=1,2,\cdots$.
Therefore
$$
\m\left(\,\bigcap_{l=1}^\infty \bigcup_{i=l}^\infty 
\left\{x:|2^i(g_i - g_{i+1})| > 1 \right\} \setminus A_{n(m)} \right)
=0,
$$
that is,
$$
\m\left(\,\bigcap_{l=1}^\infty \bigcup_{i=l}^\infty
\left\{x:|g_i - g_{i+1}| > {1 \over 2^i}\right\} 
\setminus A_{n(m)} \right) = 0.
$$

Let $L(1) = 1$, and for each $p=2,3,\cdots$ find $L(p) \geq L(p-1)$
such that 
$$
\m\left(\,\bigcup_{i=L(p)}^\infty
\left\{x:|g_i - g_{i+1}| > {1\over 2^i} \right\}
\setminus A_{n(m)} \right) \leq {1 \over p}.
$$
Define
$$
B_p = \bigcup_{i=L(p)}^\infty
\left\{x:|g_i - g_{i+1}| > {1 \over 2^i} \right\} 
\setminus A_{n(m)},\ p=1,2,\cdots . 
$$
Observe that $B_1 \supset B_2 \supset B_3 \supset \cdots,$ and
$\m(\bigcap_{p=1}^\infty B_p) = 0$.
Suppose $x \notin \bigcap_{p=1}^\infty B_p$ and
$x \notin A_{n(m)}$.
We will show that $(g_i(x))_{i=1}^\infty$ converges in this case.
First there is a $p_x$ such that $x \notin B_{p_x}$.
Therefore $|g_i(x) - g_{i+1}(x)| \leq 1/2^i$ for all
$i \geq L(p_x)$.
Let $\a > 0$ be given.  Find $M$ such that $2/2^M \leq \a$ and
$M \geq L(p_x)$. Suppose $j > i \geq M$. Then
$$
\eqalign{
|g_i(x) - g_j(x)| &\leq |g_i(x) - g_{i+1}(x)| + 
|g_{i+1}(x) - g_{i+2}(x)| + \cdots + |g_{j-1}(x) - g_j(x)| \cr 
&\leq {1 \over 2^i} + {1 \over 2^{i+1}} + \cdots + {1 \over 2^{j-1}}
\cr
&< {2 \over 2^i} \leq {2 \over 2^M} \leq \a.
}$$
So $(g_i(x))_{i=1}^\infty$ is a Cauchy sequence in $\R$.
For all $x$ 
define
$$
g_k^m(x) = \cases
\lim_{i \to \infty} g_i(x), &
 x \notin \bigcap_{p=1}^\infty B_p
            \cup A_{n(m)} \\
0, & x \in \bigcap_{p=1}^\infty B_p \cup A_{n(m)}. \\
\endcases
$$
\noindent
$g_k^m$ is the pointwise limit of measurable functions, namely
$$
g_i\cdot \chi_{\left(\bigcap_{p=1}^\infty B_p\right)^c \cap (A_{n(m)})^c},
$$
so $g_k^m$ is measurable.
Let $B_k^m = \bigcap_{p=1}^\infty B_p$.

We now remember that $k$ and $m$ were arbitrarily chosen, so for 
each $w_k^m$ we have defined $g_k^m$ and $B_k^m$ for $k=1,\cdots ,2^m; 
m=m_0,m_0+1,\cdots.$
Let 
$$
B = \bigcup_{n=1}^\infty \bigcup_{k=1}^{2^m} B_k^m.
$$
$\m(B) = 0$ since $B$ is the countable union of sets with zero
measure.
Define
$$
S(v) = \lim_{m \to \infty} \sum_{k=1}^{2^m} g_k^m.
$$
It is not immediately clear that $S(v)$ exists. 
We turn to this question next.

We claim that for almost all $x \notin A_{n(m)}$ ($m \geq m_0$) we have
$$
g_p^m(x) = g_{2p-1}^{m+1}(x) + g_{2p}^{m+1}(x).
$$
Proof of claim: We know that 
$$
\left((g_p^m)_i\right)_{i=1}^\infty \subset T(w_p^m) 
\hbox{\rm\ converges in } (A_{n(m)})^c,
$$
$$
\left((g_{2p-1}^{m+1})_i\right)_{i=1}^\infty \subset T(w_{2p-1}^{m+1}) 
\hbox{\rm\ converges in } (A_{n(m+1)})^c \supset (A_{n(m)})^c,
$$
$$
\left((g_{2p}^{m+1})_i\right)_{i=1}^\infty \subset T(w_{2p}^{m+1}) 
\hbox{\rm\ converges in } (A_{n(m+1)})^c \supset (A_{n(m)})^c, \hbox{\rm and }
$$
$$
\left(\left(g_{2p-1}^{m+1}\right)_i + \left(g_{2p}^{m+1}\right)_i\right)_{i=1}^\infty
\subset T\left(w_p^m\right)
$$
because $T$ is additive.
(The notation $((g_p^m)_i)_{i=1}^\infty$ simply means the sequence
$(g_i)_{i=1}^\infty$ that is associated with $w_p^m$.)
Since $\|4^i(g_p^m)_i\|_0 \leq (1 + 1/2i)\cdot \|4^i\cdot T(w_p^m)\|_{L_0/E}
\leq 2\cdot \sigma(T(w_p^m)) \leq 2 \d_{n(m)}/6$ we have
$$\left\|4^i\left((g_p^m)_i - (g_{2p-1}^{m+1})_i - (g_{2p}^{m+1})_i\right)\right\|_0
  \leq {\d_{n(m)} \over 3} + {\d_{n(m+1)} \over 3} + {\d_{n(m+1)} \over 3}
  \leq \d_{n(m)}.
$$
Therefore
$$
\m\left(\,\bigcap_{l=1}^\infty \bigcup_{i=l}^\infty \left\{
x:\left|(g_p^m)_i(x) - (g_{2p-1}^{m+1})_i(x) - (g_{2p}^{m+1})_i(x)  \right|
> {1 \over 2^i} \right\} \setminus A_{n(m)} \right) = 0.
$$
So for almost all $x \in (A_{n(m)})^c$ 
$$
\lim_{i \to \infty} (g_p^m)_i(x) = 
\lim_{i \to \infty} \left( (g_{2p-1}^m)_i(x) + (g_{2p}^m)_i(x) \right).
$$
Thus for almost all $x \in (A_{n(m)})^c$
$$
g_p^m(x) = g_{2p-1}^{m+1}(x) + g_{2p}^{m+1}(x),
$$
which finishes the proof of the claim.

So the sequence
$$\left(\,\sum_{k=1}^{2^m} g_k^m\right)_{m=1}^\infty
$$
remains essentially fixed in $L_0((A_{n(m)})^c)$   
for $r \geq m \geq m_0$. 
So the sequence converges in 
$L_0(\bigcup_{m=1}^\infty (A_{n(m)})^c) = L_0[0,1]$, and
$S(v)$ is well-defined.

Next we show that $T=QS$. Consider $m \geq m_0$. Then
$$
\left\|\left(\,\sum_{k=1}^{2^m} g_k^m\right) 
- S(v)\right\|_0 \leq {1 \over n(m)},
$$
since the two functions are essentially identical except 
possibly on $A_{n(m)}$ and $\m(A_{n(m)}) \leq 1/n(m)$.
For each $k$ we can find $f_k^m \in T(w_k^m)$ so that
$$
\left\|g_k^m\Bigr|_{(A_{n(m)})^c} - f_k^m\Bigr|_{(A_{n(m)})^c}\right\|_0
\leq {1\over 4^m}.
$$
We then have the following inequalities:
$$
\left\|\left(\, \sum_{k=1}^{2^m} g_k^m\Bigr|_{(A_{n(m)})^c}\right) - 
\left(\, \sum_{k=1}^{2^m} f_k^m\Bigr|_{(A_{n(m)})^c}\right)
\right\|_0 \leq 2^m\cdot {1\over 4^m} 
= {1 \over 2^m},
$$
$$
\left\|\, \sum_{k=1}^{2^m} g_k^m - 
\sum_{k=1}^{2^m} f_k^m
\right\|_0 \leq  
{1 \over 2^m} + {1 \over n(m)},
$$
$$
\left\|S(v) -
\sum_{k=1}^{2^m} f_k^m
\right\|_0 \leq {1 \over 2^m} + {2 \over n(m)}.
$$
Notice that $1/2^m + 1/n(m) \to 0$ as $m \to \infty$.
The function
$$
\sum_{k=1}^{2^m} f_k^m
$$
is an element of $T(v)$. So we
can find functions in $T(v)$ that are arbitrarily
close to $S(v)$ which means that $S(v) \in T(v)$
since $E$ is closed. That is, $QS(v) = T(v)$.

Next we will show that $S$ is a continuous linear operator. 
If $S$ is additive and continuous at
zero then $S$ must also be homogeneous, and thus linear.
So it suffices to show that $S$ is additive and
continuous at zero.

$S$ is additive. To see this let $u,v \in L_0$ and let $\a > 0$ be given.
Find $m$ so that $\m(A_{n(m)}) \leq \a$. (Recall that $\m(A_{n(m)}) \leq 
1/{n(m)}$.) We will consider $v\cdot \ind,\  
u\cdot \ind,$ and $(u + v)\cdot \ind$
for an arbitrary $k$ between $1$ and $2^m$.
>From our earlier construction we have $(f_i)_{i=1}^\infty
\subset T(u\cdot \ind)$ such that $f_i \to S(u \cdot \ind)$ on 
$(A_{n(m)})^c$ and
$$
\left\|4^i\cdot f_i\right\|_0 \leq \left(1 + {1\over 2i}\right)
\left\|4^i \cdot T(u\cdot \ind)\right\|_{L_0/E},
$$
and $(g_i)_{i=1}^\infty
\subset T(v\cdot \ind)$ such that $g_i \to S(v \cdot \ind)$ on 
$(A_{n(m)})^c$ and
$$
\left\|4^i\cdot g_i\right\|_0 \leq \left(1 + {1\over 2i}\right)
\left\|4^i \cdot T(v\cdot \ind)\right\|_{L_0/E},
$$
and $(h_i)_{i=1}^\infty
\subset T((u+v)\cdot \ind)$ such that $h_i \to S((u+v)\cdot \ind)$ on 
$(A_{n(m)})^c$ and
$$
\left\|4^i\cdot h_i\right\|_0 \leq \left(1 + {1\over 2i}\right)
\left\|4^i \cdot T((u+v)\cdot \ind)\right\|_{L_0/E}.
$$
We have $f_i + g_i \in T((u+v)\cdot \ind)$ for all $i = 1,2,\cdots$.
For $i \geq 1$,
$$
\left\|4^i(f_i + g_i) - 4^i\cdot h_i\right\|_0 \leq
\left\|4^i\cdot f_i\right\|_0 + \left\|4^i\cdot g_i\right\|_0
+ \left\|4^i \cdot h_i \right\|_0 
\leq (3 + {3 \over 2i}) \d_{n(m)}/6 < \d_{n(m)}.
$$
Therefore,
$$
\m\left(\, \bigcap_{l=1}^\infty\bigcup_{i=l}^\infty
\left\{x:\left|(f_i + g_i) - h_i\right| > 
{1 \over 2^i} \right\} \setminus A_{n(m)} \right) = 0.
$$
This implies that $(f_i + g_i)$ and $h_i$ converge to the same
function on $(A_{n(m)})^c$. Thus for all $k = 1,\cdots,2^m,\ 
S(u\cdot \ind) + S(v\cdot \ind) = S((u+v)\cdot \ind)$ on 
$(A_{n(m)})^c$. Therefore $S(u) + S(v) = S(u+v)$ on $(A_{n(m)})^c$ and
$\|S(u) + S(v) - S(u+v)\|_0 \leq \a$.
Since $\a > 0$ was arbitrary we have $S(u) + S(v) = S(u+v)$.

$S$ is continuous at zero. To see this, suppose $(v_j)_{j=1}^\infty$ 
is a sequence in $L_0$ such that $v_j \to 0$. Let $\a > 0$ be given.
Find $m$ so that $1/n(m) \leq \a$. Our set $A_{n(m)}$ then has measure less
than $\a$, and $\d_{n(m)}$ is a positive number such that the closed convex
hull of the $\d_{n(m)}$-ball in $E$ is contained in the $(\a/80)$-ball in $L_0$.
There also is an $\e_{n(m)} > 0$ so that $\|f\|_0 \leq \e_{n(m)}
\Rightarrow \|Tf\|_{L_0/E} \leq \d_{n(m)}/6$, and we have $1/2^m \leq \e_{n(m)}$.
Let $j \geq 1$ be given.
For each $k = 1, \cdots ,2^m$ there is a sequence 
$$
\left( g_{j,i}^{(k)}\right)_{i=1}^\infty \subset 
  T\left(v_j \cdot \chi_k^m\right)
$$
such that
$g_{j,i}^{(k)} \to S(v_j \cdot \chi_k^m)$ on $(A_{n(m)})^c$ as $i \to \infty$ and
$$\left\|4^i\cdot4^k\cdot g_{j,i}^{(k)}\right\|_0 \leq
\left(1 + {1 \over 2i}\right)\cdot
 \left\|4^i\cdot4^k\cdot T\left(v_j\cdot
\chi_k^m\right)\right\|_{L_0/E}.
$$ 
For each $i = 1,2,\cdots$ and $k=1,\cdots,2^m$ let 
$f_{i,k} = 2^i\cdot2^k(g_{j,i}^{(k)} - g_{j,i+1}^{(k)})$.
Then $f_{i,k} \in E$ for all $i$ and $k$ and
$$
\eqalign{
\left\|2^i\cdot2^k\cdot f_{i,k}\right\|_0
&=    \left\|4^i\cdot4^k\cdot (g_{j,i}^{(k)} - g_{j,i+1}^{(k)})\right\|_0 \cr
&\leq \left\|4^i\cdot4^k\cdot g_{j,i}^{(k)}\right\|_0 
      + \left\|4^{i+1}\cdot4^k\cdot g_{j,i+1}^{(k)}\right\|_0 \cr
&\leq 3\cdot\sigma\left(T\left(v_j\cdot\chi_k^m\right)\right) \leq \d_{n(m)}.
}$$
Using the technique employed in proving Lemma 1.2 we can conclude that
$$
\m\left(\,\bigcup_{k=1}^{2^m}\bigcup_{i=1}^\infty\left\{
   x:\left|g_{j,i}^{(k)} - g_{j,i+1}^{(k)}\right| > 
   {1 \over 2^i} \cdot {1 \over 2^k}
\right\}\right) \leq \a.
$$
Let the set above be called $D$ (so $\m(D) \leq \a$). 
Find $I$ such that $2/2^I \leq \a$. Then
$$
\left\|S\left(v_j\cdot\chi_k^m\right) - g_{j,I}^{(k)}\right\|_{L_0((A_{n(m)})^c
  \cup D^c)}
  \leq \sum_{i=I}^\infty {1 \over 2^i} \cdot {1 \over 2^k} = {2 \over 2^I} \cdot
  {1 \over 2^k}.
$$
Therefore
$$
\left\|S(v_j) - \sum_{k=1}^{2^m} g_{j,I}^{(k)}\right\|_{L_0((A_{n(m)})^c
  \cup D^c)} \leq \sum_{k=1}^{2^m} {2 \over 2^I} \cdot {1 \over 2^k}
  = {2 \over 2^I} \leq \a,
$$
and
$$
\left\|S(v_j) - \sum_{k=1}^{2^m} g_{j,I}^{(k)}\right\|_0 \leq 3\a.
$$
This is true for any $j \geq 1$.
Now
$$
\left\|\sum_{k=1}^{2^m} g_{j,I}^{(k)}\right\|_0
 \leq
 2\sum_{k=1}^{2^m} \left\|4^I\cdot 4^k \cdot 
 T\left(v_j\cdot \ind\right)\right\|_{L_0/E}.
$$
Since $T$ is continuous for each $k$, $\left\|4^I\cdot 4^k \cdot 
 T\left(v_j\cdot \ind\right)\right\|_{L_0/E}$ goes to zero as 
$j$ goes to infinity. Therefore the whole sum goes to zero as
$j$ goes to infinity. 
So $\limsup_{j \to \infty} \|S(v_j)\|_0 \leq 3 \a$. However, $\a > 0$ was
arbitrary, so $\lim_{j \to \infty} S(v_j) = 0$. 
That is, $S$ is  a continuous linear operator.

Suppose that $S'$ is another continuous linear operator from
$L_0$ to $L_0$ such that $QS' = T$.
Then $Q(S - S') = QS - QS' = T - T = 0$, whence $S - S'$ maps
$L_0$ into the locally convex
\line{space $E$.  We conclude that $S = S'$. \hfil
$\square$}

The proof of Theorem 2.1 works with a milder assumption on the subspace
$E$.  It does not have to be locally convex - the key assumption is only
that given a neighborhood $V$ of $0$ there is a smaller neighborhood $U$
so that if $x_n \in U$ then $\sum_{n=1}^N 2^{-n} x_n$ is in $V$ for all
$N$ (i.e. $E$ is exponentially galbed in the sense of Turpin \cite{8}).  We can
generalize further by replacing the sequence $(2^{-n})$ with a strictly
positive term sequence $(a_n)$ such that $\sum a_n < \infty$.  By a
classical result due to Aoki \cite{1} and Rolewicz \cite{7} 
we know that locally bounded
spaces are locally $p$-convex for some $p>0$. Also, if $U$ is locally
$p$-convex then $\sum_{n=1}^N 2^{-(n/p)} U \subset U$ for all $N$.  In this
way we can see that the generalized result includes locally bounded
subspaces of $L_0$.

We can combine Theorem 2.1 with Kwapien's theorem \cite{4}.

\bigskip

\proclaim{Theorem 2.2} Let $S: L_0 \to L_0$ be a linear operator.
Then
$$S(f)(x) = \sum_{n=1}^\infty g_n(x) f(\sigma_n(x))
$$
for every $f \in L_0$, where
\settabs 16 \columns
\+ & (i) & each $\sigma_n :[0,1] \to [0,1]$ is a non-singular 
measurable map, \cr
\+ & (ii) & each $g_n$ is in $L_0$, \cr
\+ & (iii) & for almost all $x$ in $[0,1]$, $g_n(x) \neq 0$
for only finitely many $n$. \cr

\noindent
Conversely, every map defined in the above way is a linear
operator from $L_0$ to $L_0$.
\endproclaim

\proclaim{Corollary 2.3} Let $E$ be a locally convex subspace of
$L_0$ and $Q$ be the quotient map. 
Then $T$ is an operator from $L_0$ to $L_0/E$ if
and only if $T = QS$ for some $S$ of the form in Theorem 2.2. 
\endproclaim

By following the proof of Theorem 4.1 in \cite{2} we have the 
following corollary.

\proclaim{Corollary 2.4} Let $E$ and $F$ be subspaces of $L_0$, each of which
is either locally convex or locally bounded.  Then $L_0/E$ is isomorphic
to $L_0/F$ if and only if there is an isomorphism $S$ of $L_0$ to itself
such that $S(E) = F$.
\endproclaim

\vskip2truecm

\document

\vfill
\eject

\enddocument